\documentclass[11pt]{amsart}
\usepackage[left=3cm,right=3cm,top=3cm,bottom=3cm]{geometry}
\usepackage{latexrelease}
\usepackage{enumerate, mathrsfs,color,amssymb,orcidlink}
\usepackage{amsmath,latexsym}
\usepackage[varg]{pxfonts}
\usepackage{hyperref}
\IfFileExists{fontawesome5.sty}{\usepackage{fontawesome5}}{\newcommand{\faIcon}[1]{}}
\usepackage{xcolor}
\usepackage{enumitem}
\usepackage{booktabs}

\newtheorem{theorem}{Theorem}[section]

\theoremstyle{definition}

\newtheorem{example}[theorem]{Example}
\newtheorem{corollary}[theorem]{Corollary}

\numberwithin{equation}{section}
\def\sqr#1#2{{\,\vcenter{\vbox{\hrule height.#2pt\hbox{\vrule width.#2pt
height#1pt \kern#1pt\vrule width.#2pt}\hrule height.#2pt}}\,}}

\definecolor{orcidcolor}{HTML}{A6CE39}

\usepackage[numbers,sort&compress]{natbib}
\allowdisplaybreaks[4]

\begin{document}
\title[Diagonals of Partial Isometries]{Diagonals of Partial Isometries with Prescribed Kernel and Cokernel Dimensions}

\author[Z. Wang]{Zhenquan Wang}
\address[Z. Wang]{School of Mathematical Sciences, Qufu Normal University, Qufu, Shandong, 273165, China}
\email{wangzhenquan1997@163.com}

\author[L. Wang]{Liguang Wang$^*$}
\address[L. Wang]{School of Mathematical Sciences, Qufu Normal University, Qufu, Shandong, 273165, China}
\email{wangliguang0510@163.com}

\thanks{$^*$ Corresponding author. E-mail: wangliguang0510@163.com}

\thanks{ORCID iDs:  \href{https://orcid.org/0009-0009-0976-8466}{https://orcid.org/0009-0009-0976-8466}(Zhenquan Wang),  \href{https://orcid.org/0000-0001-6449-7405}{https://orcid.org/0000-0001-6449-7405}(Liguang Wang)}
\begin{abstract}
We characterize the diagonals of partial isometries on a separable infinite-dimensional Hilbert space with prescribed kernel and cokernel dimensions. This answers Question 6.2 of Jasper, Loreaux, and Weiss.
\end{abstract}
\maketitle

\section{Introduction}

Throughout, $\mathcal{H}$ is a separable infinite-dimensional complex Hilbert space and let $\mathcal B(\mathcal{H})$ denote the spaces of bounded linear operators on $\mathcal{H}$. For $A\in\mathcal B(\mathcal{H})$, let $A^*$ denote its adjoint and put $|A|=(A^*A)^{1/2}$. For $u,v\in \mathcal{H}$, define the operator $u\otimes v$ by
$$
 (u\otimes v)x=\langle x,v\rangle u\qquad(x\in \mathcal{H}).
$$

For $A\in\mathcal B(\mathcal{H})$ and an orthonormal basis $e=(e_j)_{j\ge1}$ of $\mathcal{H}$, define
$$
 \operatorname{diag}_e A=(\langle Ae_j,e_j\rangle)_{j\ge1}.
$$
The set of diagonals of $A$ is
$$
 \mathcal D(A)=\{\operatorname{diag}_e A:e\text{ is an orthonormal basis of }\mathcal{H}\}.
$$

Fix a reference orthonormal basis $e=(e_j)_{j\ge1}$ of $\mathcal{H}$, and write $\operatorname{diag}A=\operatorname{diag}_e A$. For a bounded scalar sequence $a=(a_j)_{j\ge1}$, $\operatorname{diag}(a)$ is the operator defined by $\operatorname{diag}(a)e_j=a_je_j$. The unitary orbit of $A$ is
$$
 \mathcal U(A)=\{UAU^*:U\text{ is unitary on }\mathcal{H}\}.
$$
Thus
$$
 \mathcal D(A)
 =\{\operatorname{diag}B:B\in\mathcal U(A)\}.
$$

An operator $A\in\mathcal B(\mathcal{H})$ is positive, denoted by $A\ge0$, if $A=A^*$ and $\langle Ax,x\rangle\ge0$ for every $x\in \mathcal{H}$. For $A\ge0$, its trace is
$$
 \operatorname{Tr}A=\sum_{j\ge1}\langle Ae_j,e_j\rangle\in[0,\infty].
$$

Let $A$ be a self-adjoint $N\times N$ matrix with eigenvalues $\lambda_1\ge\cdots\ge\lambda_N$. By the Schur--Horn theorem \cite{Schur,Horn}, a real vector $d$ is a diagonal of a unitary conjugate of $A$ if and only if its nonincreasing rearrangement $d^\downarrow$ satisfies
$$
 \sum_{j=1}^n d_j^\downarrow\le\sum_{j=1}^n\lambda_j
 \quad(1\le n<N),\qquad
 \sum_{j=1}^N d_j^\downarrow=\sum_{j=1}^N\lambda_j .
$$
The matrix with this spectrum and diagonal can be chosen real symmetric by Horn's construction.

For $A\in\mathcal B(\mathcal{H})$, define its real and imaginary parts by
$$
 \operatorname{Re}A=\frac{A+A^*}{2},\qquad
 \operatorname{Im}A=\frac{A-A^*}{2i}.
$$
 
Let $a=(a_j)_{j\ge1}$ be a nonnegative sequence converging to zero. Arrange its positive terms in nonincreasing order and denote the resulting sequence by $a^\downarrow$. If $a$ has only $r<\infty$ positive terms, set $a_n^\downarrow=0$ for $n>r$. For nonincreasing sequences $d,t$ of nonnegative numbers converging to zero, we say that $d$ is majorized by $t$, and write $d\prec t$, if
\begin{equation}\label{eq:maj}
 \sum_{j=1}^n d_j\le\sum_{j=1}^n t_j\quad(n\ge1),\qquad
 \sum_{j\ge1}d_j=\sum_{j\ge1}t_j\quad\text{in }[0,\infty].
\end{equation}

We write $\|\cdot\|_2$ for the Hilbert--Schmidt norm. An operator is called real when its matrix in the specified basis has real entries. For operators between spaces with specified bases we use the same convention.

Neumann \cite{Neumann} extended the Schur--Horn convexity theorem to infinite dimensions by taking closures in the $\ell^\infty$ norm. Arveson and Kadison \cite{AK} studied diagonals of self-adjoint operators, including positive trace-class operators, and formulated related questions for von Neumann algebras. Kaftal and Weiss \cite{KW} developed majorization methods that also apply to positive compact operators outside the trace class. They characterized the diagonals of injective positive compact operators and showed that the admissible diagonals can be realized by orthogonal conjugation of a real diagonal operator.

For diagonals of orthogonal projections, Kadison \cite{KadisonI,KadisonII} showed that when the sums of the diagonal entries at most $1/2$ and of the deficits from $1$ for entries greater than $1/2$ are both finite, their difference must be an integer. Bownik and Jasper \cite{BJCarpenter} gave a constructive proof of Kadison's theorem. Kaftal and Loreaux \cite{KL} interpreted the integer as the essential codimension of a pair of projections.

For compact self-adjoint operators, Bownik and Jasper \cite{BJ} used positive and negative excesses to compare the eigenvalue sequences with a proposed diagonal. They also studied diagonals of self-adjoint operators whose essential spectrum contains at least two points \cite{BJII}.

For bounded operators, M\"uller and Tomilov \cite[Corollary~4.3]{MT} showed that a sequence in the interior of the essential numerical range of an operator is a diagonal of that operator if the sum of the distances of its entries to the complement of that range diverges. Loreaux and Weiss \cite{LWSurvey} survey these and related diagonal problems.

Thompson \cite{Thompson} characterized the diagonals of finite matrices with prescribed singular values. Jasper, Loreaux, and Weiss \cite{JLW} extended Thompson's theorem to compact operators and characterized the diagonals of unitary operators. In \cite[Question~6.2]{JLW}, they asked for a characterization of the diagonals of partial isometries with prescribed kernel and cokernel dimensions $1\le n,m\le\infty$, with particular emphasis on the case $n=m$.

An operator $V\in\mathcal B(\mathcal{H})$ is a partial isometry if $V^*V$ is an orthogonal projection. Its kernel and cokernel dimensions are $\dim\ker V$ and $\dim\ker V^*$, respectively. For a sequence $d=(d_j)_{j\ge1}$ with $|d_j|\le1$, set
$$
 \Delta(d)=\sum_{j\ge1}(1-|d_j|)\in[0,\infty].
$$
Denote
$$
 \mathcal P(d)=\{(\dim\ker V,\dim\ker V^*):V\text{ is a partial isometry with }\operatorname{diag}V=d\}.
$$
By Theorem~\ref{thm:PI},
$$
 \mathcal P(d)=
 \begin{cases}
 \{0,1,2,\ldots,\infty\}^2,
     &\Delta(d)=\infty,\\
\{(k,k):k\in\mathbb Z,\ 1\le k\le\Delta(d)\}\cup Z(d),
     &\Delta(d)<\infty,
 \end{cases}
$$
where
$$
 Z(d)=
 \begin{cases}
 \{(0,0)\},&2\sup_j(1-|d_j|)\le\Delta(d),\\
 \varnothing,&\text{otherwise}.
 \end{cases}
$$
The criterion for $(0,0)\in\mathcal P(d)$ is the unitary diagonal theorem of Jasper, Loreaux, and Weiss \cite[Theorem~4.3]{JLW}.

Theorem~\ref{thm:PI} answers Question~6.2 of Jasper, Loreaux, and Weiss \cite{JLW}, which concerns the case $1\le n,m\le\infty$. It also covers the cases in which exactly one of $n$ and $m$ is zero. The case $n=m=0$ is the unitary diagonal theorem \cite[Theorem~4.3]{JLW}.
\section{Diagonals of partial isometries}\label{sec:PI}

Let $\mathbb D=\{z\in\mathbb C:|z|<1\}$. For a sequence $d=(d_j)_{j\ge1}$ with $d_j\in\overline{\mathbb D}$ for every $j\ge1$, set
$$
 \alpha_j=1-|d_j|,\qquad
 \Delta(d)=\sum_{j=1}^{\infty}\alpha_j\in[0,\infty].
$$

For $T\in\mathcal B(\mathcal{H})$, the essential numerical range is
$$
 W_e(T)=\left\{\lambda\in\mathbb C:
 \langle Tx_r,x_r\rangle\to\lambda
 \text{ for some orthonormal sequence }(x_r)_{r\ge1}\subset \mathcal{H}
 \right\}.
$$
The set $W_e(T)$ is compact and convex; see \cite[Section~3]{MT}.

We use the real finite-rank Schur--Horn theorem \cite[Lemma~6.3 and Proposition~6.4]{KW} (see also \cite[Theorem~2.6]{JLW}). Let $r\ge1$ be an integer and let $\beta=(\beta_j)_{j\ge1}$ satisfy $\beta_1\ge\cdots\ge\beta_r>0$ and $\beta_j=0$ for $j>r$. For a nonnegative summable sequence $\alpha=(\alpha_j)_{j\ge1}$, there exists a positive operator $B$ of rank $r$ on $\mathcal{H}$ with real matrix entries in the fixed basis $(e_j)_{j\ge1}$, positive eigenvalues $\beta_1,\ldots,\beta_r$ counted with multiplicity, and diagonal entries $\langle Be_j,e_j\rangle=\alpha_j$ for every $j\ge1$ if and only if $\alpha^\downarrow\prec\beta$ in the sense of \eqref{eq:maj}.

\begin{theorem}\label{thm:PI}
Let $d=(d_j)_{j\ge1}$ be a complex sequence, and let $n,m\in\{0,1,2,\ldots,\infty\}$. There exists a partial isometry $V\in\mathcal B(\mathcal{H})$ with diagonal $d$ and
$$
 \dim\ker V=n,\qquad \dim\ker V^*=m
$$
if and only if $|d_j|\le1$ for every $j\ge1$ and one of the following conditions holds:
\begin{enumerate}[label=\textup{(\roman*)},leftmargin=2.1em]
\item $\Delta(d)=\infty$.
\item $\Delta(d)<\infty$ and $1\le n=m\le\Delta(d)$.
\item $\Delta(d)<\infty$, $n=m=0$, and
\begin{equation}\label{eq:unitarycriterion}
 2\sup_{j\ge1}(1-|d_j|)\le\Delta(d).
\end{equation}
\end{enumerate}
\end{theorem}

\begin{proof}
Let $(e_j)_{j\ge1}$ be the fixed orthonormal basis of $\mathcal{H}$. Every partial isometry $V$ is a contraction, so $|\langle Ve_j,e_j\rangle|\le1$ for every $j$. Thus we may assume $|d_j|\le1$ for every $j\ge1$. Let $a_j=|d_j|$ and $\alpha_j=1-a_j$ for $j\ge1$. Set $a=(a_j)_{j\ge1}$, $\alpha=(\alpha_j)_{j\ge1}$, and $\Delta=\Delta(d)$.

Define the diagonal unitary $D=\operatorname{diag}(\omega_j)$ by
$$
 \omega_j=
 \begin{cases}
 \overline{d_j}/|d_j|,&d_j\ne0,\\
 1,&d_j=0.
 \end{cases}
$$
For a partial isometry $V$ with diagonal $d$,
$$
 \operatorname{diag}(DV)=(\omega_jd_j)_{j\ge1}=a,\qquad
 (DV)^*(DV)=V^*V.
$$
Thus $DV$ is a partial isometry. Since $D$ is unitary,
$$
 \ker(DV)=\ker V,\qquad
 \ker((DV)^*)=D(\ker V^*),
$$
so both kernel dimensions are unchanged. Conversely, if $W$ is a partial isometry with diagonal $a$, then
$$
 \operatorname{diag}(D^*W)
 =(\overline{\omega_j}a_j)_{j\ge1}=d,
$$
and left multiplication by $D^*$ also preserves both kernel dimensions. Since the conditions on $d$ in the theorem depend only on $|d_j|$, it suffices to prove the theorem for $a$.

We first prove necessity. Suppose that $V$ is a partial isometry with diagonal $a$ and
$$
 \dim\ker V=n,\qquad \dim\ker V^*=m.
$$
Since $\Delta=\infty$ is condition~\textup{(i)}, we only need to consider $\Delta<\infty$. Let $P_0$ and $Q_0$ be the orthogonal projections onto $\ker V$ and $\ker V^*$, respectively. Since $V^*V=I-P_0$,
\begin{equation}\label{eq:HSidentity}
 \begin{aligned}
 \|(V-I)e_j\|^2
 &=\|Ve_j\|^2+1-2\operatorname{Re}\langle Ve_j,e_j\rangle\\
 &=\langle(I-P_0)e_j,e_j\rangle+1-2a_j\\
 &=2\alpha_j-\langle P_0e_j,e_j\rangle
 \le2\alpha_j.
 \end{aligned}
\end{equation}
Hence $R=V-I$ is Hilbert--Schmidt, with
$$
 \|R\|_2^2
 =\sum_{j=1}^{\infty}\|(V-I)e_j\|^2
 \le2\Delta<\infty.
$$
Since $V=I+R$,
$$
 P_0=I-V^*V=-R-R^*-R^*R,\qquad
 Q_0=I-VV^*=-R-R^*-RR^*.
$$
Since $R$ is compact, $P_0$ and $Q_0$ are compact orthogonal projections and hence have finite rank. Moreover,
$$ 
 n-m
 =\operatorname{Tr}(P_0-Q_0)
 =\operatorname{Tr}(RR^*)-\operatorname{Tr}(R^*R)=0. 
$$
Thus $n=m<\infty$.

Set $C=I-\operatorname{Re}V$. Since $\|V\|\le1$,
$$
 \langle Cx,x\rangle
 =\|x\|^2-\operatorname{Re}\langle Vx,x\rangle
 \ge\|x\|^2-\|Vx\|\,\|x\|
 \ge0
 \qquad(x\in \mathcal{H}).
$$
Thus $C\ge0$. Since $\langle Ce_j,e_j\rangle=1-a_j=\alpha_j$,
$$
 \sum_{j=1}^{\infty}\|C^{1/2}e_j\|^2
 =\sum_{j=1}^{\infty}\langle Ce_j,e_j\rangle
 =\sum_{j=1}^{\infty}\alpha_j
 =\Delta<\infty.
$$
Hence $C^{1/2}$ is Hilbert--Schmidt, $C$ is trace class, and $\operatorname{Tr}C=\Delta$.

Since $P_0$ is the orthogonal projection onto $\ker V$, we have $VP_0=0=P_0V^*$. Therefore,
$$
 P_0CP_0
 =P_0-\tfrac12(P_0VP_0+P_0V^*P_0)
 =P_0.
$$
Since $0\le P_0\le I$, we have $0\le C^{1/2}P_0C^{1/2}\le C$. Thus
\begin{equation}\label{eq:lowerdefect}
 n=\operatorname{Tr}P_0
  =\operatorname{Tr}(P_0CP_0)
  =\operatorname{Tr}(C^{1/2}P_0C^{1/2})
  \le\operatorname{Tr}C=\Delta.
\end{equation}
If $n\ge1$, then \textup{(ii)} holds. If $n=0$, then $m=0$ and $V^*V=VV^*=I$, so $V$ is unitary. By \cite[Theorem~4.3]{JLW},
$$
 2\sup_{j\ge1}(1-a_j)\le\Delta.
$$
Thus \textup{(iii)} holds. This completes the proof of necessity.

We now prove sufficiency. Suppose that one of \textup{(i)}--\textup{(iii)} holds. If \textup{(iii)} holds, then by \cite[Theorem~4.3]{JLW} there exists a unitary operator $W$ with diagonal $a$. The operator $D^*W$ is unitary and has diagonal $d$, so both kernel dimensions are zero.

For \textup{(ii)}, set $k=n=m$. Then $1\le k\le\Delta<\infty$. Set
$$
 \Delta-k=2\ell+\eta,\qquad
 \ell\in\{0,1,2,\ldots\},\qquad
 0\le\eta<2,
$$
and
$$
 c=1-\eta,\qquad
 x=(1,0)^{\mathsf T},\qquad
 y=(c,\sqrt{1-c^2})^{\mathsf T}.
$$
Let $\mathcal{H}_0=\mathbb C^2\oplus\mathbb C^{k-1}\oplus\mathbb C^\ell \oplus\ell^2$, define $A$ by
\begin{equation}\label{eq:PImodel}
 A=(x\otimes y)\oplus0_{\mathbb C^{k-1}}
   \oplus(-I_{\mathbb C^\ell})\oplus I_{\ell^2}.
\end{equation} 
The summand $0_{\mathbb C^{k-1}}$ is omitted when $k=1$, and $-I_{\mathbb C^\ell}$ is omitted when $\ell=0$. Since $\|x\|=\|y\|=1$,
$$
 (x\otimes y)^*(x\otimes y)=y\otimes y,\qquad
 (x\otimes y)(x\otimes y)^*=x\otimes x.
$$
Thus $x\otimes y$ is a partial isometry on $\mathbb C^2$ with one-dimensional kernel and cokernel. The other summands are also partial isometries. Hence $A$ is a real partial isometry and
$$
 \dim\ker A=\dim\ker A^*=1+(k-1)=k.
$$

Set $B=I_{\mathcal{H}_0}-\operatorname{Re}A$. For a real unitary $Q:\mathcal{H}_0\to \mathcal{H}$, the operator $QAQ^*$ is real and
$$
 \operatorname{diag}(QAQ^*)
 =\operatorname{diag}\bigl(\operatorname{Re}(QAQ^*)\bigr)
 =\operatorname{diag}(I_{\mathcal{H}}-QBQ^*).
$$
Since $a_j=1-\alpha_j$ for every $j\ge1$, it suffices to find a real unitary $Q:\mathcal{H}_0\to \mathcal{H}$ with $\operatorname{diag}(QBQ^*)=\alpha$.

With respect to the decomposition of $\mathcal{H}_0$,
$$
 B=B_2\oplus I_{\mathbb C^{k-1}}
       \oplus2I_{\mathbb C^\ell}\oplus0_{\ell^2},
 \qquad
 B_2=
 \begin{pmatrix}
 1-c&-\sqrt{1-c^2}/2\\
 -\sqrt{1-c^2}/2&1
 \end{pmatrix}.
$$
Since $c=1-\eta$,
$$
 \det(zI_2-B_2)
 =\left(z-\frac{\eta}{2}\right)
  \left(z-1-\frac{\eta}{2}\right).
$$
The eigenvalues of $B_2$ are $\eta/2$ and $1+\eta/2$. Therefore $B$ is positive and has finite rank. Its positive eigenvalues are the nonzero terms of the following nonincreasing list:
\begin{equation}\label{eq:PIbeta}
 \underbrace{2,\ldots,2}_{\ell},\quad
 1+\eta/2,\quad
 \underbrace{1,\ldots,1}_{k-1},\quad
 \eta/2.
\end{equation}
In particular,
$$
 \operatorname{Tr}B
 =2\ell+1+\frac{\eta}{2}+k-1+\frac{\eta}{2}
 =k+2\ell+\eta
 =\Delta.
$$

Denote these positive eigenvalues by $\beta_1\ge\cdots\ge\beta_r>0$ and set $\beta_j=0$ for $j>r$. Since $\sum_{j=1}^\infty\alpha_j=\Delta<\infty$, we have $\alpha_j\to0$. Let $\alpha^\downarrow$ be the positive entries of $\alpha$ arranged in nonincreasing order, with zeros appended if there are only finitely many positive entries. Since $0\le\alpha_j\le1$ and $\beta_j\ge1$ for $j<r$,
$$
 \sum_{j=1}^s\alpha_j^\downarrow
 \le s\le\sum_{j=1}^s\beta_j
 \qquad(1\le s<r).
$$
For $s\ge r$,
$$
 \sum_{j=1}^s\alpha_j^\downarrow
 \le\Delta
 =\sum_{j=1}^r\beta_j
 =\sum_{j=1}^s\beta_j.
$$
Moreover, $\sum_{j=1}^{\infty}\alpha_j^\downarrow=\Delta=\sum_{j=1}^{\infty}\beta_j$. Hence $\alpha^\downarrow\prec\beta$.

By the real finite-rank Schur--Horn theorem, there exists a real positive finite-rank operator $B_1$ on $\mathcal{H}$ with diagonal $\alpha$ and positive eigenvalues $\beta$.  Since $B$ and $B_1$ have finite rank on infinite-dimensional spaces, both kernels are infinite-dimensional. Choose real orthonormal bases $(g_j)_{j\ge1}$ of $\mathcal{H}_0$ and $(h_j)_{j\ge1}$ of $\mathcal{H}$ such that
$$
 Bg_j=\beta_jg_j,\qquad
 B_1h_j=\beta_jh_j
 \qquad(j\ge1),
$$
where $\beta_j=0$ for $j>r$. Define $Q:\mathcal{H}_0\to \mathcal{H}$ by $Qg_j=h_j$ for every $j\ge1$. Since $Q$ maps a real orthonormal basis onto a real orthonormal basis, it is real and unitary. Moreover,
$$
 B_1Qg_j=\beta_jh_j=QBg_j
 \qquad(j\ge1).
$$
Thus $B_1Q=QB$, and
$$
 B_1=QBQ^*,\qquad
 \operatorname{diag}(QBQ^*)=\alpha.
$$
Thus
$$
 \operatorname{diag}(QAQ^*)
 =\operatorname{diag}(I_{\mathcal{H}}-B_1)
 =(1-\alpha_j)_{j\ge1}
 =a.
$$
Set $V=D^*QAQ^*$. Since $D$ and $Q$ are unitary, $V$ is a partial isometry and
$$
 \operatorname{diag}V
 =(\overline{\omega_j}a_j)_{j\ge1}=d,\qquad
 \dim\ker V=\dim\ker V^*=k.
$$
This proves sufficiency in \textup{(ii)}.

It remains to prove sufficiency in \textup{(i)}. Suppose that $\Delta=\infty$. Let $S$ be the unilateral forward shift on $\ell^2(\mathbb N)$, let $U$ be the bilateral shift on $\ell^2(\mathbb Z)$, and define
$$
 V_0=U\oplus S^{(m)}\oplus(S^*)^{(n)}.
$$
Here $S^{(m)}$ and $(S^*)^{(n)}$ denote the orthogonal direct sums of $m$ copies of $S$ and $n$ copies of $S^*$, respectively. The corresponding summand is omitted when $m=0$ or $n=0$, and $\infty$ denotes countably many copies. Since $S$ and $S^*$ are partial isometries and $U$ is unitary, $V_0$ is a partial isometry. Moreover,
$$
 \ker S=\{0\},\qquad
 \dim\ker S^*=1,\qquad
 \ker U=\ker U^*=\{0\}.
$$
Hence
$$
 \dim\ker V_0=n,\qquad
 \dim\ker V_0^*=m.
$$

We next verify that $W_e(V_0)=\overline{\mathbb D}$. Let $(u_j)_{j\in\mathbb Z}$ be the standard orthonormal basis of the first summand $\ell^2(\mathbb Z)$, so that $Uu_j=u_{j+1}$ for every $j\in\mathbb Z$. Fix $\lambda\in\mathbb C$ with $|\lambda|=1$ and define
$$
 \xi_r
 =2^{-r/2}\sum_{j=2^r}^{2^{r+1}-1}\lambda^{-j}u_j
 \qquad(r\ge1).
$$
For $r<s$, we have $2^{r+1}-1<2^s$. Hence $\langle u_j,u_q\rangle=0$ whenever $2^r\le j\le2^{r+1}-1$ and $2^s\le q\le2^{s+1}-1$. Therefore,
$$
 \langle\xi_r,\xi_s\rangle
 =2^{-(r+s)/2}
   \sum_{j=2^r}^{2^{r+1}-1}
   \sum_{q=2^s}^{2^{s+1}-1}
   \lambda^{-j}\overline{\lambda^{-q}}
   \langle u_j,u_q\rangle=0.
$$
Since $\langle\xi_s,\xi_r\rangle=\overline{\langle\xi_r,\xi_s\rangle}$, we obtain $\langle\xi_r,\xi_s\rangle=0$ for every $r\ne s$. Moreover, since $|\lambda|=1$,
$$
 \|\xi_r\|^2
 =2^{-r}\sum_{j=2^r}^{2^{r+1}-1}|\lambda^{-j}|^2
 =2^{-r}\cdot2^r
 =1.
$$
Thus $(\xi_r)_{r\ge1}$ is an orthonormal sequence. We also have
$$
 \langle V_0\xi_r,\xi_r\rangle
 =\langle U\xi_r,\xi_r\rangle
 =2^{-r}\sum_{j=2^r+1}^{2^{r+1}-1}
   \lambda^{-(j-1)}\overline{\lambda^{-j}}
 =2^{-r}(2^r-1)\lambda
 =(1-2^{-r})\lambda
 \longrightarrow\lambda
$$
Thus $\lambda\in W_e(V_0)$. Since $\lambda$ was arbitrary with $|\lambda|=1$, every point on the unit circle belongs to $W_e(V_0)$. Since $W_e(V_0)$ is convex,
$$
 \overline{\mathbb D}
 =\operatorname{conv}\{z\in\mathbb C:|z|=1\}
 \subseteq W_e(V_0).
$$
Moreover, for every unit vector $x$,
$$
 |\langle V_0x,x\rangle|
 \le\|V_0x\|
 \le\|V_0\|=1.
$$
Thus $W_e(V_0)\subseteq\overline{\mathbb D}$. Hence $W_e(V_0)=\overline{\mathbb D}$.

Let
$$
 J=\{j\ge1:a_j<1\},\qquad
 \mathcal{H}_J=\overline{\operatorname{span}}\{e_j:j\in J\}.
$$
Since $W_e(V_0)=\overline{\mathbb D}\subseteq\overline{W(V_0)}$, where $W(V_0)$ is the numerical range of $V_0$, the real affine hull of $W(V_0)$ is $\mathbb C$. For every $j\in J$, we have $a_j\in\operatorname{int}W_e(V_0)=\mathbb D$ and $\operatorname{dist}(a_j,\mathbb C\setminus W_e(V_0))=1-a_j$. Since $a_j=1$ for $j\notin J$, we obtain
$$
 \sum_{j\in J}\operatorname{dist}\bigl(a_j,\mathbb C\setminus W_e(V_0)\bigr)
 =\sum_{j\in J}(1-a_j)
 =\Delta=\infty.
$$
Each term is at most $1$, so $J$ is infinite. Enumerate $J$ as $j_1<j_2<\cdots$. By \cite[Corollary~4.3]{MT}, there is an orthonormal basis $(f_r)_{r\ge1}$ of the space on which $V_0$ acts such that $\langle V_0f_r,f_r\rangle=a_{j_r}$ for every $r\ge1$.

Let $Q$ be the unitary onto $\mathcal{H}_J$ defined by $Qf_r=e_{j_r}$, and set $T_J=QV_0Q^*$. Then $T_J$ is a partial isometry on $\mathcal{H}_J$, unitarily equivalent to $V_0$, and
$$
 \langle T_Je_{j_r},e_{j_r}\rangle
 =\langle V_0f_r,f_r\rangle
 =a_{j_r}\qquad(r\ge1).
$$
On $\mathcal{H}=\mathcal{H}_J\oplus \mathcal{H}_J^\perp$, define
$$
 W=T_J\oplus I_{\mathcal{H}_J^\perp}.
$$
Then $W$ is a partial isometry. Since $T_J$ has diagonal $(a_j)_{j\in J}$ and $a_j=1$ for $j\notin J$, we have $\operatorname{diag}W=a$. Also
$$
 \dim\ker W=\dim\ker T_J=n,\qquad
 \dim\ker W^*=\dim\ker T_J^*=m.
$$
Set $V=D^*W$. Since $D$ is unitary,
$$
 V^*V=W^*W,\qquad
 \ker V=\ker W,\qquad
 \ker V^*=D^*(\ker W^*).
$$
Thus $V$ is a partial isometry, and
\begin{equation*}
 \operatorname{diag}V
 =(\overline{\omega_j}a_j)_{j\ge1}=d,\qquad
 \dim\ker V=n,\qquad
 \dim\ker V^*=m.
 \qedhere
\end{equation*}
\end{proof}

\begin{corollary}\label{cor:PIequality}
Let $V$ be a partial isometry with diagonal $d$ and suppose that
$$
 \dim\ker V=\dim\ker V^*=k,\qquad 1\le k<\infty,\qquad \Delta(d)=k.
$$
For every diagonal unitary $D$ with $\operatorname{diag}(DV)=(|d_j|)_{j\ge1}$, the operator $DV$ is an orthogonal projection with kernel dimension $k$.
\end{corollary}

\begin{proof}
Fix such a diagonal unitary $D$ and set $V_1=DV$. Let $P$ be the orthogonal projection onto $\ker V_1$, and set $C=I-\operatorname{Re}V_1$. Since $D$ is unitary, $\ker V_1=\ker V$, so $\operatorname{rank}P=k$.

Since $V_1$ is a contraction, $C=I-\operatorname{Re}V_1\ge0$. Its diagonal is $(1-|d_j|)_{j\ge1}$, so $C$ is trace class and $\operatorname{Tr}C=\Delta(d)=k$. Moreover, $V_1P=0$ and $PV_1^*=0$, so $PCP=P$. Hence $\operatorname{Tr}(PCP)=\operatorname{Tr}P=k$. Thus
$$
 0=\operatorname{Tr}C-\operatorname{Tr}(PCP)
  =\operatorname{Tr}\bigl((I-P)C(I-P)\bigr)
  =\|C^{1/2}(I-P)\|_2^2.
$$
It follows that $C^{1/2}(I-P)=0$. Hence $C(I-P)=(I-P)C=0$, so $C=PCP=P$.

For $x\in(\ker V_1)^\perp$, we have $\|V_1x\|=\|x\|$ because $V_1$ is a partial isometry. Thus
$$
 \|V_1x-x\|^2
 =2\|x\|^2-2\operatorname{Re}\langle V_1x,x\rangle
 =2\langle Cx,x\rangle=0.
$$
Since $V_1$ vanishes on $P\mathcal{H}$ and acts as the identity on $(I-P)\mathcal{H}$, we conclude that $DV=V_1=I-P$.
\end{proof}

The following example illustrates Theorem~\ref{thm:PI} for three diagonal sequences. The first sequence also illustrates the equality case in Corollary~\ref{cor:PIequality}.

\begin{example}\label{ex:PI}
Recall that $\mathcal P(d)$ is the set of pairs $(\dim\ker V,\dim\ker V^*)$ for partial isometries $V$ with diagonal $d$.

\begin{enumerate}[label=\textup{(\arabic*)},leftmargin=2.1em]
\item Let $d_j=1-2^{-j}$ for $j\ge1$. Then
$$
 \Delta(d)=\sum_{j=1}^{\infty}2^{-j}=1,\qquad
 2\sup_{j\ge1}(1-|d_j|)=2\cdot\frac12=1.
$$
Condition~\textup{(ii)} of Theorem~\ref{thm:PI} holds exactly when $n=m=1$. Condition~\textup{(iii)} holds when $n=m=0$, with equality in~\eqref{eq:unitarycriterion}. Hence
$$
 \mathcal P(d)=\{(0,0),(1,1)\}.
$$

Let $V$ be any partial isometry with diagonal $d$ and $\dim\ker V=1$. Then $\dim\ker V^*=1$. Moreover,
$$
 d=|d|,\qquad
 \Delta(d)=\dim\ker V=\dim\ker V^*=1.
$$
By Corollary~\ref{cor:PIequality} with $D=I$, $V$ is an orthogonal projection. Thus $d$ is the diagonal of both a unitary operator and an orthogonal projection with one-dimensional kernel.

\item Let $d=(1/4,1,1,\ldots)$. Then
$$
 \Delta(d)=\frac34,\qquad
 2\sup_{j\ge1}(1-|d_j|)=\frac32.
$$
Since $\Delta(d)<1$, condition~\textup{(ii)} of Theorem~\ref{thm:PI} cannot hold. The inequality in condition~\textup{(iii)} also fails, since $3/2>3/4$. Therefore,
$$
 \mathcal P(d)=\varnothing.
$$
This sequence has all its terms in the closed unit disk, but it is not the diagonal of any partial isometry.

\item Let $d_j=1/2$ for every $j\ge1$. Then
$$
 \Delta(d)=\sum_{j=1}^{\infty}\frac12=\infty.
$$
By part~\textup{(i)} of Theorem~\ref{thm:PI}, for every $n,m\in\{0,1,2,\ldots,\infty\}$ there exists a partial isometry $V$ with diagonal $d$ and
$$
 \dim\ker V=n,\qquad
 \dim\ker V^*=m.
$$
Consequently,
$$
 \mathcal P(d)=\{0,1,2,\ldots,\infty\}^2.
$$
For this diagonal, the kernel and cokernel dimensions can be prescribed independently, including unequal dimensions.
\end{enumerate}
\end{example}
\noindent\textbf{Funding.} This work is partially supported by the National Natural Science Foundation of China (12471124) and the Natural Science Foundation of Shandong Province, China (ZR2025MS93).\\

\noindent\textbf{Author contributions.} All authors performed the research and prepared the first draft of the manuscript. All authors read and approved the final manuscript.\\

\noindent\textbf{Data availability.} This article has no associated data.\\

\noindent\textbf{Competing interests.} The authors declare that they have no conflict of interest.\\


\begin{thebibliography}{99}

\bibitem{AK}
W. Arveson and R. V. Kadison,
\emph{Diagonals of self-adjoint operators},
in Operator Theory, Operator Algebras, and Applications,
Contemp. Math. \textbf{414}, American Mathematical Society,
Providence, RI, 2006, 247--263.

\bibitem{BJCarpenter}
M. Bownik and J. Jasper,
\emph{Constructive proof of the Carpenter's theorem},
Canad. Math. Bull. \textbf{57} (2014), 463--476.

\bibitem{BJ}
M. Bownik and J. Jasper,
\emph{Diagonals of self-adjoint operators I: Compact operators},
J. Funct. Anal. \textbf{289} (2025), no.~5, Paper No.~110939.

\bibitem{BJII}
M. Bownik and J. Jasper,
\emph{Diagonals of self-adjoint operators II: non-compact operators},
Math. Ann. \textbf{391} (2025), 431--507.

\bibitem{Horn}
A. Horn,
\emph{Doubly stochastic matrices and the diagonal of a rotation matrix},
Amer. J. Math. \textbf{76} (1954), 620--630.

\bibitem{JLW}
J. Jasper, J. Loreaux, and G. Weiss,
\emph{Thompson's theorem for compact operators and diagonals of unitary operators},
Indiana Univ. Math. J. \textbf{67} (2018), no.~1, 1--27.

\bibitem{KadisonI}
R. V. Kadison,
\emph{The Pythagorean theorem. I. The finite case},
Proc. Natl. Acad. Sci. USA \textbf{99} (2002), 4178--4184.

\bibitem{KadisonII}
R. V. Kadison,
\emph{The Pythagorean theorem. II. The infinite discrete case},
Proc. Natl. Acad. Sci. USA \textbf{99} (2002), 5217--5222.

\bibitem{KL}
V. Kaftal and J. Loreaux,
\emph{Kadison's Pythagorean theorem and essential codimension},
Integral Equations Operator Theory \textbf{87} (2017), 565--580.

\bibitem{KW}
V. Kaftal and G. Weiss,
\emph{An infinite dimensional Schur--Horn theorem and majorization theory},
J. Funct. Anal. \textbf{259} (2010), 3115--3162.

\bibitem{LWSurvey}
J. Loreaux and G. Weiss,
\emph{On diagonals of operators: selfadjoint, normal and other classes},
in \#operatortheory27: Conference Proceedings, Timi\c{s}oara,
July 2--6, 2018, Theta Foundation, Bucharest, 2020, 193--214.

\bibitem{MT}
V. M\"uller and Y. Tomilov,
\emph{Diagonals of operators and Blaschke's enigma},
Trans. Amer. Math. Soc. \textbf{372} (2019), 3565--3595.

\bibitem{Neumann}
A. Neumann,
\emph{An infinite-dimensional version of the Schur--Horn convexity theorem},
J. Funct. Anal. \textbf{161} (1999), 418--451.

\bibitem{Schur}
I. Schur,
\emph{\"Uber eine Klasse von Mittelbildungen mit Anwendungen auf die Determinantentheorie},
Sitzungsber. Berl. Math. Ges. \textbf{22} (1923), 9--20.

\bibitem{Thompson}
R. C. Thompson,
\emph{Singular values, diagonal elements, and convexity},
SIAM J. Appl. Math. \textbf{32} (1977), 39--63.

\end{thebibliography}
\end{document}